\documentclass[12pt,leqno]{amsart}   
\usepackage{amsmath}
\usepackage{amssymb}

\newcommand{\bea}[1]{\begin{eqnarray}\label{#1}}
\newcommand{\eea}{\end{eqnarray}}
\newcommand{\beas}{\begin{eqnarray*}}
\newcommand{\eeas}{\end{eqnarray*}}

\newcounter{rocount}

\newcounter{gracount}

\newtheorem{theor}{Theorem}[section]
\newtheorem{propo}[theor]{Proposition}
\newtheorem{corol}[theor]{Corollary}
\newtheorem{lemma}[theor]{Lemma}

\theoremstyle{definition}

\theoremstyle{remark}

\renewcommand{\Box}{\blacksquare}

\newcommand{\NN}{{\mathbb{N}}}
\newcommand{\ZZ}{{\mathbb{Z}}}

\renewcommand{\epsilon}{\varepsilon}
\renewcommand{\phi}{\varphi}
\newcommand{\norm}[1]{{\left\|{#1}\right\|}}

\newcommand{\cstar}{\mbox{$C^*$}}

\newcommand{\id}{\operatorname{id}}

\newcommand{\BBB}{{\mathbb B}}

\newcommand{\TT}{{\mathbb T}}
\newcommand{\FF}{{\mathbb F}}
\newcommand{\Log}{{\operatorname{Log}}}
\newcommand{\myA}{A_{\epsilon}}
\newcommand{\myE}{E_{\alpha}}

\newcommand{\mydual}{{{\alpha}}}
\newcommand{\myB}{B_{\epsilon}}
\newcommand{\matrM}{{\mathbf M}}
\newcommand{\HHH}{{\mathcal H}}
\newcommand{\interU}{{T_{0,m}}}
\newcommand{\interH}{{K_{n,m}}}

\newcommand{\NEW}[1]{}   
\reversemarginpar

\begin{document}
\title{ Finite dimensional representations of the soft torus}
\author{S\o ren Eilers}
\thanks{Partially supported by the Carlsberg Foundation}
\address{Matematisk Afdeling\\
K\o benhavns Universitet\\
Universitetsparken 5\\
DK-2100 Copenhagen \O\\
Denmark}
\email{eilers@math.ku.dk}
\urladdr{http://www.math.ku.dk/\~{}eilers}
\author{Ruy Exel}
\address{
 Departamento de Matematica\\
   Universidade Federal de Santa Catarina\\
   Campus Universitario - Trindade\\
   88010-970 Florianopolis SC\\
   Brazil}
\email{exel@mtm.ufsc.br}
\urladdr{http://www.mtm.ufsc.br/\~{}exel}
\date{\today}
\begin{abstract}
The soft tori constitute a continuous deformation, in a very precise
sense, from the commutative \cstar-algebra $C(\TT^2)$ to the highly
non-commutative \cstar-algebra $C^*(\FF_2)$. Since both of these
\cstar-algebras are known to have a separating family of finite dimensional
representations, it is natural to ask whether that is also the case
for the soft tori. We show that this is in fact the case.
\end{abstract}
\maketitle

\section{Introduction}
Knowing that a given \cstar-algebra has many 
representations on finite dimensional Hilbert spaces is of great
importance to understanding structural properties of it.
Among
\cstar-algebras, those who possess a separating family of
finite dimensional  representations are called \emph{residually
finite dimensional} or just \emph{RFD}.  This class was studied in
\cite{krgpm:frfdc}, \cite{retal:frfpc} and \cite{rja:rfc}, and more
recent insight about it 
has lead to important advances in  classification theory and the
theory of quasidiagonal \cstar-algebras (see, e.g., \cite{bbek:gilfdc},
\cite{md:nsa} and \cite{md:aqc}).

For any $\epsilon\geq 0$ we define a \cstar-algebra $\myA$ as the
unital universal \cstar-algebra defined by the generators $u,v$ subject to 
the relations
\begin{gather*}
uu^*=u^*u=1\qquad vv^*=v^*v=1\qquad \norm{uv-vu}\leq \epsilon.
\end{gather*}
As recorded in \cite{re:staacm},  $A_0$ is  the commutative \cstar-algebra of functions
over the torus $\TT^2$, and $\myA$ is the full
\cstar-algebra of the free group of two generators $\FF_2$ whenever
$\epsilon\geq 2$. For
$\epsilon$ between $0$ and $2$ we get a class of \cstar-algebras which
are commonly referred to as \emph{soft tori}.
These \cstar-algebras are of relevance to several problems in
operator algebra theory (see \cite{retal:iacu}) and have been
extensively studied in \cite{cc:ndck}, 
 \cite{re:staacm}, \cite{re:stvacn},
\cite{gaeretal:stiii}. 

The starting point of the investigation reported on in the present
paper is a result from
\cite{re:stvacn}, stating that the soft tori form a continuous field
interpolating between the (hard) torus and the group \cstar-algebra
of the free group.
Since $C(\TT^2)$ is obviously RFD, and since $\cstar(\FF_2)$ was
proved to be RFD in \cite{mdc:fcfgtg} --- a surprise at the time --- we are naturally led to the
question of whether the same is true for the interpolating family
$\myA$. We are going  to prove that this is the case.

\section{Methods}
We are going to employ a method from \cite{cc:ndck} and \cite{re:staacm} which lies
behind many results about the structural properties of $\myA$. 

We define $\myB$ as the universal \cstar-algebra given by the generators
$\{u_n\}_{n\in\ZZ}$ and the relations
\begin{gather}\label{myBrels}
u_nu^*_n=u^*_nu_n=1\qquad \norm{u_{n+1}-u_{n}}\leq \epsilon.
\end{gather}
Clearly one can define an automorphism $\mydual$ on $\myB$ by 
\[
u_n\mapsto u_{n+1},
\]
and
as seen in
\cite{re:staacm}  one has
\[
\myA=\myB\rtimes_\mydual\ZZ.
\]
There is a faithful conditional expectation $\myE:\myA\longrightarrow \myB$.

Our strategy is to
prove that $\myA$ is RFD by proving that $\myB$ is RFD in a way which
is covariant with $\mydual$.

\section{ Finite dimensional representations of $\myB$}

We start out by finding a new picture of $\myB$ by generators and
relations.

\begin{lemma}\label{retire}
For any $\epsilon<2$, $\myB$ is isomorphic to the universal
\cstar-algebra generated by $v_0, \{h_n\}_{n\in \ZZ}$ subject to the
relations
\begin{gather}\label{myBaltrelations}
v_0v_0^*=v_0^*v_0=1\qquad
h_n=h_n^*\qquad
\|h_n\|\leq 2 \cos(\varepsilon/2)
\end{gather}
\end{lemma}
\begin{proof}
Let us denote the \cstar-algebra generated by $v_0$ and $h_n$ subject
to \eqref{myBaltrelations}
by $\myB'$. 
We can define a map $\phi:\myB\longrightarrow\myB'$  by
\begin{gather*}
u_n\mapsto\left\{\begin{array}{ll}
  e^{i\pi  h_n}\cdots e^{i\pi
h_1}v_0&n>0\\
v_0&n=0\\
  e^{-i\pi  h_n}\cdots e^{-i\pi
h_{-1}}v_0&n<0	 \end{array}\right.
\end{gather*}
since the elements to the right of the arrow above 
satisfy the relations \eqref{myBrels}.
Similarly, 
the universal property of $\myB'$ allows for a map
$\psi:\myB'\longrightarrow \myB$ defined by
\begin{gather*}
v_0\mapsto u_0\qquad h_n\mapsto \dfrac1{i\pi}\Log(u_{n+1}u_{n}^*)
\end{gather*}
Clearly $\phi$ and $\psi$ are each others' inverse.
\end{proof}

This characterization can be used to shorten the proof of
\cite[2.2]{re:staacm}, stating that $\myB$ is homotopic to
$C(\TT)$. To see this, define maps $\phi:\myB'\longrightarrow C(\TT)$
and $\psi:C(\TT)\longrightarrow\myB'$ by the correspondence
$v_0\leftrightarrow [z\mapsto z]$, $h_n\leftrightarrow 0$. Clearly $\phi\psi=\id_{C(\TT)}$, and
$\chi_t:\myB'\longrightarrow\myB'$ given by
\[
\chi_t(v_0)=v_0\qquad\chi_t(h_n)=t h_n
\]
provides a homotopy from $\id_{\myB'}$ to $\psi\phi$.

In the following proof, we denote by  $\operatorname{Alg}(X)$ 
the smallest $*$-algebra, not necessarily  closed, generated by $X$ inside some \cstar-algebra.

\begin{propo}\label{BisRFD}
For any $\epsilon<2$, $\myB$ is RFD. In fact, for any $0\not=b\in \myB$ there
exists $n\in\NN$, an automorphism $\beta$ of $\matrM_n$  and a representation $\rho:\myB\longrightarrow
\matrM_n$ with the properties:
\[
\rho(b)\not=0 \qquad
\beta\rho=\rho\mydual\qquad
\beta^n=\id.
\]
\end{propo}
\begin{proof}
For the first claim we use the characterization of $\myB$
given by the relations \eqref{myBaltrelations} and proceed as in
\cite{mdc:fcfgtg}. Fix a faithful non-degenerate representation $\pi:\myB\longrightarrow \BBB(\HHH)$, where we
may assume that $\HHH$ is a separable Hilbert space.

Let $P_m$ be a 
sequence of projections, with $\operatorname{rank}(P_m)=m$, converging strongly to the unit of
$\BBB(\HHH)$, and abbreviate
\[
\interU=P_m\pi(v_0)P_m\qquad\interH=P_m\pi(h_n)P_m
\]
Now note that for each $m$ the collection of elements
$\{V_{0,m},H_{n,m}: n\in \ZZ\}$  defined by 
\[
V_{0,m}=  \left[\begin{matrix} \interU & \sqrt{P_m -
\interU\interU^*}\\ \sqrt{P_m - 
\interU^*\interU} & \interU^*\end{matrix}\right]\quad
H_{n,m}=\left[  \begin{matrix} \interH&0\\0&\interH\end{matrix}\right]
\]
satisfies \eqref{myBaltrelations} in{\NEW{$\matrM_2$!}} $\matrM_2(P_m\BBB(\HHH) P_m)\simeq \matrM_{2m}$. Consequently we get
 representations $\pi_m:\myB\longrightarrow\matrM_{2m}$. We are going
to check, following \cite{mdc:fcfgtg}, that
\[
\underline{\pi}:\myB\longrightarrow\prod_{m=1}^\infty{\matrM_{2m}}\qquad
\underline{\pi}(b)=({\pi_m(b)})_{m=1}^\infty
\]
is an isometry. It suffices to check that
$\norm{\underline{\pi}(x)}\geq \norm{x}-\eta$ for any $\eta>0$ and any 
$x\in\operatorname{Alg}(\{v_0,h_{-N},\dots,h_N\})$. Fix $\eta$, $N$
and $x$ and  write
\[
x=F(v_0,h_{-N},\dots,h_N)
\]
where $F$ is some finite linear combination of finite words in $2N+2$
variables and their adjoints. Since
\[
V_{0,m}\longrightarrow\left[\begin{matrix}\pi(v_0)&0\\0&\pi(v_0)^*\end{matrix}\right]\qquad
H_{n,m}\longrightarrow\left[\begin{matrix}\pi(h_n)&0\\0&\pi(h_n)\end{matrix}\right],\qquad
m\longrightarrow \infty
\]
strongly on the unit ball of{\NEW{$\matrM_2$!}} $\matrM_2(\BBB(\HHH))$ we conclude that
\begin{eqnarray*}
\lefteqn{\lim_{m\longrightarrow\infty}\norm{F(V_{0,m},H_{-N,m},\dots,H_{N,m})}}\\
&&=\norm{
\left[\begin{matrix}\pi(F(v_0,h_{-N},\dots,h_{N}))&0\\
0&\pi(F(v_0^*,h_{-N},\dots,h_{N}))\end{matrix}\right]}\\
&&\geq
\norm{\pi(x)}=\norm{x}.
\end{eqnarray*}
We can hence find $m$ such that
\[
\norm{\underline{\pi}(x)}\geq
\norm{{\pi_m}(x)}=\norm{F(V_{0,m},H_{-N,m},\dots,H_{N,m})}\geq \norm{x}-\eta.
\]

For the second claim, we go back to the original
presentation \eqref{myBrels} of $\myB$ by unitary generators only. For a given
$b\in\myB$ with $\norm{b}=1$ we fix, using the first part of the proposition, a finite dimensional representation
$\pi:\myB\longrightarrow \matrM_{m}$ with $\norm{\pi(b)}>\tfrac34$. We also
fix $c$ and $N\in\NN$ such that
\[
\norm{b-c}<\tfrac14\qquad c\in
 \operatorname{Alg}(u_{-N},\dots,u_{N}).
\]
Choose $M>0$ and unitaries $v^\pm_0,\dots,v^\pm_{M}$ with the
properties
\[
\norm{v^\pm_{n+1}-v^\pm_{n}}\leq \epsilon\qquad
v^\pm_0=u_{\pm N}\qquad
v^\pm_M=1.
\]
There is then exactly one representation $\pi':\myB\longrightarrow
\matrM_m$ which is $2(N+M)$-periodic in the sense that $\pi'(u_n)=\pi'(u_{n+2(N+M)})$ and satisfies
\[
\pi'(u_n)=\left\{\begin{array}{ll}v^-_{-N-n}&-M-N\leq n<-N\\
\pi(u_n)&-N\leq  n\leq N\\
v^+_{n-N}&N<n\leq N+M
\end{array}\right.
\]
Note that $\pi'(c)=\pi(c)$; in particular $\norm{\pi'(c)}\geq\tfrac12$.

Now let $n=2(N+M)m$. With $\beta$ defined  as the forward cyclic shift
in block form (with period $2(N+M)$) we  may define a covariant representation $\rho$ of $\myB$ on
$\matrM_{n}$ by
\[
  u_i \mapsto 
\left[  \begin{matrix}\pi'(u_i) \\
           & \pi'(u_{i+1}) \\
           && \ddots \\
           &&& \pi'(u_{i+2(N+M)-1}) \\
  \end{matrix}\right].
\]
We have
\[
\norm{\rho(b)}\geq \norm{\pi'(b)}\geq \norm{\pi'(c)}-\tfrac14>0
\]
\end{proof}

\section{ Finite dimensional representations of $\myA$}

\begin{theor}\label{AisRFD}
For any $\epsilon>0$, let $\myA$ be the universal \cstar-algebra generated by
a pair of unitaries subject to the relation
$\norm{uv-vu}\leq\epsilon$. Then $\myA$ is residually finite
dimensional in the sense that it admits a separating family of finite
dimensional representations. 
\end{theor}
\begin{proof}
We may assume that $0<\epsilon <2$. Let $0 \neq a\in \myA$.  Then
also $b= \myE(a^*a)$
is nonzero, for the conditional expectation is faithful. Choose $n$, $\rho$ and $\beta$ as in Proposition \ref{BisRFD} and
define
\[
\pi:\myA=\myB\rtimes_\mydual\ZZ\longrightarrow
\matrM_{n}\rtimes_\beta\ZZ
\]
as the extension to the crossed product of the covariant
$*$-homomorphism $\rho$. We then have, with $E_\beta$ the conditional
expectation from $\matrM_n\rtimes_\beta\ZZ$ to $\matrM_n$
\[
E_\beta(\pi(a^*a))=\pi(\myE(a^*a))=\rho(b)\not=0,
\]
so $\pi(a)\not=0$. 

Note finally that since $\beta$ is periodic, any irreducible
representation of $\matrM_{n}\rtimes_\beta
\ZZ$  is finite dimensional, and hence this \cstar-algebra is RFD. Therefore 
we may compose $\pi$ with some finite dimensional representation of
it to exhibit a finite dimensional representation which does not
vanish on $a$.
\end{proof}

Linear algebra tells us that $\matrM_n$ has a faithful tracial state
and has the property that every matrix which is \emph{hyponormal} is the sense that 
\[
x^*x\geq xx^*
\]
is in fact normal. As in  \cite{mdc:fcfgtg}, we may conclude:

\begin{corol}\label{choiripoff}
For any $\epsilon$, $\myA$ has a faithful tracial state. 
Furthermore, any hyponormal operator in $\myA$ is in fact
normal.
\end{corol}
\begin{proof}
These two properties clearly pass from matrices to sums of the form
$\bigoplus_{n\in\NN}{\matrM_{m_n}}$, and from these sums to any of
their subalgebras. By the
theorem, $\myA$ is one
such.
\end{proof}

\providecommand{\bysame}{\leavevmode\hbox to3em{\hrulefill}\thinspace}

\end{document}